\documentclass[11pt]{article}
\usepackage{ulem}
\usepackage[makeroom]{cancel}
\usepackage{color}
\usepackage{xcolor}
\usepackage{amsfonts}
\usepackage{amssymb}
\usepackage{amscd}
\usepackage{amsmath}
\usepackage{epsfig}
\usepackage{graphicx, subfigure}
\usepackage{latexsym}
\usepackage{pst-3dplot}
\usepackage{verbatim}
\usepackage{tikz}
\usepackage{pdfsync}
\usepackage{hyperref}
\usepackage{multirow}
\usepackage[utf8]{inputenc}

\usepackage{color}
\usepackage{graphicx}

\usepackage{caption}
\usepackage{lipsum}
\newlength{\dinwidth}
\newlength{\dinmargin}
\newtheorem{definition}{Definition}
\newtheorem{theorem}{Theorem}

\newtheorem{proposition}{Proposition}
\newtheorem{corollary}{Corollary}
\newtheorem{remark}{Remark}
\newtheorem{lemma}{Lemma}
\newtheorem{example}{Example}

\def \surf{{\cal L}(x)}

\def \i{{\rm i}}

\def\px1{p_{x_1}}
\def\px2{p_{x_2}}
\def\pu1{p_{u_1}}

\def\y{y_0}

\def\X{{\cal X}}

\def\aa{A}
\def\bb{B}

\def\u{{\hat u}}
\def\w{{\hat w}}

\makeatletter
\@addtoreset{equation}{section}
\makeatother

\begin{document}
\title{Isoperiodic deformations of Abelian differentials of the second kind over elliptic curves  and the Boussinesq equation}

\author{Vladimir Dragovi\'c$^1$ and Vasilisa Shramchenko$^2$}
\date{}

\maketitle

\footnotetext[1]{Department of Mathematical Sciences, University
	of Texas at Dallas, 800 West Campbell Road, Richardson TX 75080,
	USA. Mathematical Institute SANU, Kneza Mihaila 36, 11000
	Belgrade, Serbia.  E-mail: {\tt
		Vladimir.Dragovic@utdallas.edu}}

\footnotetext[2]{Department of mathematics, University of
	Sherbrooke, 2500, boul. de l'Universit\'e,  J1K 2R1 Sherbrooke, Quebec, Canada. E-mail: {\tt Vasilisa.Shramchenko@Usherbrooke.ca}}

\begin{abstract}
We study deformations of a genus one Riemann surface and of a second order Abelian differential on the surface which preserve the periods of the differential with respect to a chosen canonical homology basis of the surface.  We call these deformations {\it isoperiodic}. We derive a second order ordinary differential equation with rational coefficients governing the  variations of the position of the  unique pole of the differential under  the isoperiodic deformations. The obtained equation depends on the order of the pole of the differential. We characterize the solutions of the obtained ordinary differential equations that correspond to the isoperiodic deformations. We apply these results to the theory of genus one solutions to the Boussinesq equation.
\end{abstract}

\smallskip


\smallskip


\tableofcontents

\newpage

\section{Introduction}

The Painlev\'e VI equation PVI$(\alpha, \beta, \gamma, \delta)$ is a family of second order ordinary differential equations
with rational coefficients parameterized by four values $\alpha, \beta, \gamma, \delta \in \mathbb C$ \cite{Ince}. The celebrated Picard solution of the PVI$(0;0;0;\frac{1}{2})$  is obtained as follows \cite{Picard}. For two arbitrary constants $c_1, c_2\in\mathbb C$, the inverse Abel-Jacobi image of the point with the coordinates $(c_1,c_2)$
 in the Jacobian of the elliptic curve
\begin{equation}
\label{curve}
v^2=u(u-1)(u-x)
\end{equation}
gives a point $Q_0$ on the curve. Its $u$-coordinate as a function of $x$ is the Picard solution of the Painlev\'e equation with parameters $\alpha=\beta=\gamma=0$ and $\delta=\frac{1}{2}$.

In \cite{DS2019} we interpret the Picard solution as {\it isoperiodic} deformations of an Abelian differential of the third kind over a family of elliptic curves \eqref{curve} parameterized by $x\in\mathbb C\setminus\{0,1\}$, that is deformations preserving the periods of the differential.

Here we extend this idea to Abelian differentials of the second kind on elliptic curves.
Namely, we consider a meromorphic differential having a unique pole of order $n+2$ for any integer $n\geq 0$ at a given point $Q_0$ on the compact Riemann surface corresponding to the elliptic curve \eqref{curve}.
We study the deformations the curve and of the $u$-coordinate $y_0$ of the pole $Q_0$ as a function of $x$ preserving the periods of the differential with respect to a chosen canonical homology basis of the surface.  We call such deformations {\it isoperiodic}.

As a result, we obtain a second order ordinary differential equation \eqref{ode} for $y_0(x)$ for every given $n.$ Remarkably, coefficients of these equations are rational functions of $x$ and $y_0$, the property shared with the Painlev\'e equations.
To derive these equations, we use Bell polynomials \cite{parts} and Rauch variational formulas \cite{Fay92} on Riemann surfaces of genus one.  We characterize the solutions of the obtained ordinary differential equations that correspond to the isoperiodic deformations.

A similar question was studied in \cite{GS} although the equations and their consequences derived there  are very different from ours.

We apply our results to the study of genus one solutions \cite{Dub1981} of Boussinesq equation \cite{Bous}.

\section{Elliptic curves and Abelian differentials}
\label{sect_diff}

For the elliptic curve defined by equation \eqref{curve}
denote by  $\surf$  the corresponding compact Riemann surface. Let us denote an arbitrary point  of $\surf$  by $P=(u(P), v(P))$. The covering $u:\surf\to\mathbb CP^1$ sending a point $P$ to $u(P)$ is ramified at the four ramification points $P_0=(0,0)$, $P_1=(1,0)$, $P_x=(x,0)$ and $P_\infty=(\infty, \infty)$ and their projection on the $u$-sphere, $0,1,x, \infty$ are called the branch points of the covering $u$.

Introduce the {\it standard local coordinates} on the surface $\surf$ as follows:
\begin{align}
\label{coordinates}
&\zeta_{u_k}(P)=\sqrt{u(P) - u_k} \quad\mbox{if}\quad P\sim P_{u_k}, \quad\mbox{with}\quad u_k\in\{0,1,x\}\,,
\nonumber
\\
& \zeta_\infty(P)= \frac{1}{\sqrt{u(P)}}\quad\mbox{if}\quad P\sim P_\infty,
\\
& \zeta_Q(P)=u(P)-u(Q) \quad\mbox{if}\quad P\sim Q \quad\mbox{and $Q$ is a regular point.}
\nonumber
\end{align}
It is with respect to these standard local coordinates that we define the {\it evaluation} of Abelian differentials at a point of the Riemann surface. For an Abelian differential $\Upsilon$ we define its {\it value} at a point $\tilde Q\in\surf$ which is not a pole of $\Upsilon$ as the constant term of the Taylor series expansion of the differential with respect to the standard local parameter $\zeta$ from the list \eqref{coordinates} at $\tilde Q$, that is
\begin{equation}
\label{evaluation}
\Upsilon(\tilde Q) = \frac{\Upsilon(P)}{d\zeta(P)}\Big{|}_{P=\tilde Q}.
\end{equation}

Let us fix a canonical homology basis $\{a,\,b\}$ on $\surf$ and consider the holomorphic normalized differential $\omega$ on $\surf$. For this differential we can write
\begin{equation}
\label{omega}
\omega=\frac{du}{I_0\, v},
\end{equation}
where the normalization constant $I_0$ is defined by
\begin{equation}
\label{I0}
I_0=\oint_a\frac{du}{v}\,.
\end{equation}

Evaluating $\omega$ at ramification points and at some regular point $Q_0$ according to \eqref{evaluation}, yields
\begin{eqnarray}
&&\omega(P_0) := \frac{\omega(Q)}{d\sqrt{u(Q)}}{\Big |}_{Q=P_0, u=0} = \frac{2}{I_0\sqrt{x}}\,,
\label{omega_atP0}
\\
&&\omega(P_1) := \frac{\omega(Q)}{d\sqrt{u(Q)-1}}{\Big |}_{Q=P_1,u=1} = \frac{2}{I_0\sqrt{1-x}}\,,
\label{omega_atP1}
\\
&&\omega(P_x) := \frac{\omega(Q)}{d\sqrt{u(Q)-x}}{\Big |}_{Q=P_x, u=x} = \frac{2}{I_0\sqrt{x(x-1)}},
\label{omega_atPx}
\end{eqnarray}
\begin{equation}
\label{omega_atQ0}
\omega(Q_0):= \frac{\omega(Q)}{du(Q)}{\Big |}_{Q=Q_0} =  \frac{1}{I_0\sqrt{\y(\y-1)(\y-x)}}, \qquad\text{where }\quad y_0=u(Q_0).
\end{equation}
As a corollary of these expressions, we obtain the following useful fact about the holomorphic normalized differential:
\begin{equation}
\label{zero}
\omega^2(P_0)+\omega^2(P_1)+\omega^2(P_x)=0\,.
\end{equation}

Let us now introduce the fundamental Riemann bidifferential $W(P, Q)$ with $P, Q\in\surf$,  see \cite{Fay92}.  A bidifferential behaves as a differential with respect to each of the arguments $P$ and $Q$. As an example, one can think of $\omega_1(P)\omega_2(Q)$ where $\omega_j$ are two Abelian differentials on $\surf$, however a bidifferential is not always represented as such a product.  The Riemann bidifferential is the unique bidifferential on $\surf$ having the following three properties
\begin{itemize}
\item Symmetry: $W(P,Q) = W(Q,P);$
\item The only singularity being a second order pole along the diagonal $P=Q$ with the following local expansion in terms of a local parameter $\zeta$ near $P=Q$:
\begin{equation}
\label{W}
W(P,Q) \underset{P\sim Q}{=} \left( \frac{1}{(\zeta(P) - \zeta(Q))^2}  + {\cal O}(1) \right)d\zeta(P) d\zeta(Q);
\end{equation}
\item Normalization by vanishing of the $a$-period: $\oint_{a} W(P,Q) = 0.$
\end{itemize}
Due to the symmetry,  the above integral can be computed with respect to either $P$ or $Q$.
Clearly, $W$ depends on the choice of a canonical homology basis. As a consequence of this definition we have:
\begin{equation}
\label{bW}
\oint_{b}W(P,Q) = 2\pi{\rm i}\,\omega(P).
\end{equation}

We evaluate the bidifferential at a precise point with respect to one of its arguments in the same way as we do for differentials \eqref{evaluation} while fixing one of the arguments, say $P$, of $W$ and considering the bidifferential as a differential in the second argument $Q$. Thus, for a point $R\in \surf$ we obtain an Abelian differential of the second kind
\begin{equation*}
W(P, \tilde Q) = \frac{W(P, Q)}{d\zeta(Q)}\Big{|}_{Q=\tilde Q}.
\end{equation*}
with vanishing $a$-period and  a unique pole of order two at $P=\tilde Q$ with  the principal part $(\zeta^{-2}(P)+\mathcal O(1))d\zeta$ for $P\sim \tilde Q$ with respect to  the standard local parameter $\zeta$ at $\tilde Q$ from the list \eqref{coordinates}.

For the differentials of the second kind obtained in this way from $W$,  in the case when the pole is located
at one of the ramification points of $\surf\,,$ we can write somewhat explicit expressions in terms of the coordinates $u$ and $v$ from \eqref{curve}. To obtain such expressions we use the fact that the corresponding differential $W(P,P_{u_k})$ is the unique Abelian differential of the second kind with a second order pole at $P=P_{u_k}$ with  the principal part $(\zeta^{-2}_{u_k}(P)+\mathcal O(1))d\zeta_{u_k}$ for $P\sim P_{u_k}$ with respect to the standard parameter $\zeta_{u_k}$  \eqref{coordinates} and vanishing $a$-period. We need to introduce some normalization constants, similarly to $I_0$ in \eqref{omega}. More precisely, we have

\begin{eqnarray}
\label{WPP0}
&W(P,P_0)=\left( \frac{1}{\omega(P_0)u(P)}+I^0\right)\omega(P), \qquad &I^0=\oint_a\frac{\omega(P)}{\omega(P_0)u(P)};
\\
\label{WPP1}
&W(P,P_1)=\left( \frac{1}{\omega(P_1)(u(P)-1)}+I^1\right)\omega(P), \qquad &I^1=\oint_a\frac{\omega(P)}{\omega(P_0)(u(P)-1)};
\\
\label{WPPx}
&W(P,P_x)=\left( \frac{1}{\omega(P_x)(u(P)-x)}+I^x\right)\omega(P), \qquad &I^x=\oint_a\frac{\omega(P)}{\omega(P_0)(u(P)-x)}.
\end{eqnarray}
Evaluating with respect to the second argument at another ramification point,  we get
\begin{eqnarray}
\label{W0x}
&&W(P_x,P_0)=\left( \frac{1}{\omega(P_0)x}+I^0\right)\omega(P_x)
=\left( I^x-\frac{1}{\omega(P_x)x}\right)\omega(P_0),
\\
\label{W1x}
&&W(P_x,P_1)=\left( \frac{1}{\omega(P_1)(x-1)}+I^1\right)\omega(P_x)=\left( I^x-\frac{1}{\omega(P_x)(x-1)}\right)\omega(P_1)\,;
\end{eqnarray}
and, similarly,
\begin{eqnarray}
\label{WQ0x}
&& W(Q_0,P_x)=\left( \frac{1}{\omega(P_x)(\y-x)}+I^x\right)\omega(Q_0)\,.
\end{eqnarray}
Comparing the two expressions in each of \eqref{W0x} and \eqref{W1x}, we obtain
\begin{eqnarray}
\label{I0x}
&&\frac{I^0}{\omega(P_0)} = \frac{I^x}{\omega(P_x)} - \frac{xI_0^2}{4};
\\
\label{I1x}
&&\frac{I^1}{\omega(P_1)}= \frac{I^x}{\omega(P_x)} - \frac{(x-1)I_0^2}{4};
\\
\label{I01}
&&\frac{I^0}{\omega(P_0)}=\frac{I^1}{\omega(P_1)}- \frac{I_0^2}{4}.
\end{eqnarray}

\section{Isoperiodic deformations}

Let us fix an arbitrary point $Q_0\in\surf$ different from the ramification points and let $\y:=u(Q_0)$ be its projection on the $u$-sphere.
Let $\aa\in\mathbb C$ and $n\in \mathbb N\cup\{0\}$ be also fixed. Define the following Abelian differential of the second kind
\begin{equation}
\label{Omega}
\Omega(P)=W^{(n)}(P, Q_0)+\aa \omega(P)\,,\quad P\in \surf,
\end{equation}
where $\omega$ the holomorphic normalized differential on $\surf$ and
where $W^{(n)}$ stands for the Abelian differential of the second kind with zero $a$-period and the following principal  part with respect to $u$ at its only pole at $P=Q_0:$
\begin{equation}
\label{pole}
W^{(n)}(P, Q_0)\underset{P\sim Q_0}{=}\left( \frac{1}{(u(P)-\y)^{n+2}} + \mathcal O (1) \right)du(P)\,.
\end{equation}
Here $u(P)-\y$ being the local parameter near $Q_0$ and $\y=u(Q_0)\,.$ Note that for $n=0$ the differential $W^{(0)}(P, Q_0)=W(P, Q_0).$ For positive $n$, differential \eqref{pole} may be interpreted as the $n$th coefficient in the Taylor series expansion of $W(P,Q)/d\zeta_{Q_0}(Q)$ with respect to the local parameter $\zeta_{Q_0}(Q)$ \eqref{coordinates} at $Q=Q_0$ and divided by $(n+1).$

The $b$-period of the normalized meromorphic differential \eqref{pole} is given by (see, for example, Corollary 10-5 in \cite{Sp1957}):
\begin{equation}
\label{bn}
\oint_{b}W^{(n)}(P, Q_0)= \frac{2\pi \i}{(n+1)} \omega^{(n)}(Q_0)\,.
\end{equation}
Here $ \omega^{(n)}(Q_0)$ is the $n$th coefficient in the Taylor series expansion of $\omega(P)/d\zeta_{Q_0}(P)$ at $P=Q_0$
\begin{equation}
\label{omega-Taylor}
\omega(P) \underset{P\sim Q_0}{=} \left( \omega(Q_0)+\omega^{(1)}(Q_0 )\zeta_{Q_0}(P)+\dots +\omega^{(n)}(Q_0 )\zeta^n_{Q_0}(P) +\dots \right)d\zeta_{Q_0}(P)
\end{equation}
with respect to the standard local parameter  $\zeta_{Q_0}(Q)$ \eqref{coordinates} at $Q_0\,.$ The $b$-period of $W^{(n)}$ can be also obtained from \eqref{bW}.

Denote by $\tau$ the $b$-period of $\surf$, that is
\begin{equation}\label{eq:tau}
\tau=\oint_b\omega\,.
\end{equation}
 The periods of $\Omega$ \eqref{Omega} are
\begin{equation}
\label{periods}
\oint_a\Omega=\aa\,,\qquad \oint_b\Omega=\frac{2\pi\i}{n+1} \omega^{(n)}(Q_0)+\aa\tau=:\bb\,.
\end{equation}
Let us now consider a family of compact Riemann surfaces $\surf$ corresponding to the curves \eqref{curve} parameterized by the values of
$x$ varying in  some small, simply connected, open set $\X\subset\mathbb C\setminus\{0,1\}$. We assume that the variations are small enough for the canonical homology basis $\{a,b\}$ to be chosen consistently for all the surfaces, that is we assume that the projections on the $u$-sphere $u(a)$ and $u(b)$ of the homology basis are independent of $x\in \X$.

The $a$-period \eqref{periods} of $\Omega$ is equal to $A$ and is thus constant by construction. The $b$-period $\bb$ in \eqref{periods} is in general a function of $x$ and $\y\,.$ On the other hand, one can define a function $\y(x)$ by requiring that the $b$-period $\bb$ of $\Omega$ be constant.

Let us define a section $Q_0(x)$ of the family of surfaces $\surf\,,\,\, x\in \X$ by requiring that the periods of differential $\Omega$ \eqref{Omega} be independent of $x\,.$ Recall that we denote $\y=u(Q_0).$ Our goal is to obtain a differential equation for the function $\y(x)$.

\begin{definition}
\label{def-iso}
Let $\surf$ be a family of Riemann surfaces of genus one parameterized by $\;x\in\X\subset\mathbb C$ with a consistent choice of canonical homology bases $a,b$. Let $\Omega=\Omega(x)$ is a differential of the second kind with a unique pole on $\surf.$
We call a family of pairs $(\surf, \Omega), \;x\in\X,$ defined by conditions $\oint_a\Omega= {\rm const}$ and  $\oint_b\Omega= {\rm const}$  isoperiodic.  Alternatively, we can say that such a family gives an isoperiodic deformation of the pair $(\mathcal L, \Omega)$ with $\mathcal L\in \surf$.
\end{definition}

\section{Bell polynomials}
\label{sect_Bell}
Here we introduce polynomials $L$, which will help us differentiate expression \eqref{omega_atQ0} with respect to $\y,$ see also \cite{isoharmonic}. Let us define
\begin{equation}
\label{L}
L_l(z_1, \dots, z_l) =l! \sum_{\substack{p_1+2p_2+\dots +lp_l=l \\ p_1,\dots, p_l\geq 0}} \frac{(-1)^{\sum_{k=1}^l p_k}\, z_1^{p_1} z_2^{p_2}\dots z_l^{p_l}}{2^{\sum_{k=1}^l p_k} \prod_{k=1}^l p_k! \prod_{k=1}^l k^{p_k}}\,.
\end{equation}
For example, we have
\begin{equation*}
L_0=1, \qquad L_1(z_1) = -\frac{z_1}{2}, \qquad L_2(z_1, z_2) = \frac{z_1^2}{4} - \frac{z_2}{2},
\qquad L_3(z_1, z_2, z_3) = -\frac{z_1^3}{8} + \frac{3}{4}z_1z_2 -z_3\,.
\end{equation*}
The polynomials $L_l$  are related to the {\it complete exponential Bell polynomials}, see \cite{parts}, by a simple change of variables.
We will evaluate the polynomials at the following arguments given by sums over all the branch points of the curves $\surf$ \eqref{curve} for integer values of $k\geq 0\,:$
\begin{equation}
\label{Sigmak}
\Sigma_k:= \frac{1}{(-\y)^k}+\frac{1}{(1-\y)^k}+\frac{1}{(x-\y)^k}\,.
\end{equation}
One sees, for example, that $\frac{\partial}{\partial \y} \Sigma_k=k\Sigma_{k+1}.$  The following result shows an application of polynomials $L_l$ in our context.
\begin{lemma}
\label{lemma_Ll}
Let the polynomials $L_l$ be as above and the value $\omega(Q_0)$ be defined by \eqref{omega_atQ0}. Then for any integer $n$ and any integer $l\geq 0$
\begin{equation}
\label{Ll}
\frac{\partial^{l}\omega(Q_0)}{\partial \y^l} = \omega(Q_0) L_l(-\Sigma_1, \dots, -\Sigma_l)\,.
\end{equation}
\end{lemma}
{\it Proof.}
As is easy to see, the first derivative is $\frac{\partial\omega(Q_0)}{\partial \y} = \frac{1}{2}\omega(Q_0) \Sigma_1\,.$ Thus the derivatives in question can be written in the form of $\omega(Q_0)$ multiplied by some polynomial in $\Sigma_1, \dots, \Sigma_l$ of the following form
\begin{equation*}
\frac{\partial^{l}\omega(Q_0)}{\partial \y^l} = \frac{\partial^{l}}{\partial \y^l} \left\{\frac{1}{I_0\sqrt{\y(\y-1)(\y-x)}}\right\} = \omega(Q_0) \sum_{\substack{p_1+2p_2+\dots +lp_l=l \\ p_1,\dots, p_l\geq 0}} C_{p_1 p_2\dots p_l} \Sigma_1^{p_1}\Sigma_2^{p_2}\dots \Sigma_l^{p_l}
\end{equation*}
with some coefficients $C_{p_1 p_2\dots p_l}\,.$  Note that this notation assumes that the indices of a coefficient $C$ which are equal to zero and are placed at the end of the string of indices can be omitted, that is, for example,  $C_1=C_{1,0,0,0}\,.$ By examining the derivatives of $\omega(Q_0)$, one sees the following recursion relation for coefficients $C_{p_1 p_2\dots p_l}$
\begin{equation*}
C_{p_1 p_2\dots p_l} = \frac{1}{2}C_{p_1-1, p_2\dots p_{l-1}} + \sum_{n=1}^{l-1} n(p_n+1) C_{p_1\dots p_{n-1}, p_n+1,p_{n+1}-1, p_{n+2}\dots p_{l-1}}\,.
\end{equation*}
It is straightforward to verify that the coefficients
\begin{equation}
\label{Coefficients}
C_{p_1 p_2\dots p_l} = \frac{l!\, }{2^{\sum_{k=1}^l p_k} \prod_{k=1}^l p_k! \prod_{k=1}^l k^{p_k}}
\end{equation}
satisfy this recursion with the initial value $C_0 = 1\,.$
$\Box$
To shorten the expressions, we will write
\begin{equation}
\label{Ll-notation}
L_l:=L_l(-\Sigma_1, \dots, -\Sigma_l)= \sum_{\substack{p_1+2p_2+\dots +lp_l=l \\ p_1,\dots, p_l\geq 0}}\frac{ l!\,\Sigma_1^{p_1} \Sigma_2^{p_2}\dots \Sigma_l^{p_l}}{2^{\sum_{k=1}^l p_k} \prod_{k=1}^l p_k! \prod_{k=1}^l k^{p_k}}\,.
\end{equation}
\begin{corollary}
\label{cor_recursion}
The functions $L_n$ \eqref{Ll-notation} satisfy the following recursion relation (similarly to the Bell polynomials):
\begin{equation}
\label{recursion}
L_{n+1}=\sum _{k=0}^{n}\frac{n!}{2(n-k)!}L_{n-k}\Sigma_{k+1},
\end{equation}
where $\Sigma_j$ is defined by \eqref{Sigmak}.
\end{corollary}
{\it Proof.} From Lemma \ref{lemma_Ll} we have
\begin{equation*}
L_{n+1}=\frac{1}{\omega(Q_0)}\frac{\partial^{n+1}\omega(Q_0)}{\partial \y^{n+1}}  = \frac{1}{\omega(Q_0)}\frac{\partial^{n}}{\partial \y^{n}}\frac{\partial\omega(Q_0)}{\partial \y}
= \frac{1}{\omega(Q_0)}\frac{\partial^{n}\left\{\omega(Q_0)L_1\right\}}{\partial \y^{n}}\,.
\end{equation*}
Using the Leibniz rule for differentiating the product, we rewrite the above relation as
\begin{equation*}
L_{n+1}= \frac{1}{\omega(Q_0)} \sum_{k=0}^n {n\choose k} \frac{\partial^{n-k}\omega(Q_0)}{\partial \y^{n-k}} \frac{\partial^{k}L_1}{\partial \y^{k}}\,.
\end{equation*}
 Using now   $L_1=\frac{1}{2}\Sigma_1$ according to \eqref{Ll-notation} and  $\frac{\partial^k}{\partial \y^k}\Sigma_1=k!\Sigma_{k+1},$ which follows from definition \eqref{Sigmak} of $\Sigma_k$ and Lemma \ref{lemma_Ll} again, we prove the corollary.
$\Box$
\begin{corollary}  For the polynomials $L_n$  defined by \eqref{Ll-notation} with arguments $\Sigma_k$ given by  \eqref{Sigmak}, we have
\begin{equation}
\label{Ly0}
\frac{\partial}{\partial \y} L_n=L_{n+1}-L_nL_1\,.
\end{equation}
\end{corollary}
{\it Proof.}  Lemma \ref{lemma_Ll} implies
\begin{equation*}
\frac{\partial}{\partial \y} L_n=\frac{\partial}{\partial \y}\left\{ \frac{1}{ \omega(Q_0)}\frac{\partial^{n}\omega(Q_0)}{\partial \y^n}  \right\}
=\frac{1}{\omega(Q_0)}\frac{\partial^{n+1}\omega(Q_0)}{\partial \y^{n+1}} - \frac{1}{\omega^2(Q_0)}\frac{\partial^{n}\omega(Q_0)}{\partial \y^{n}}\frac{\partial\omega(Q_0)}{\partial \y}
\end{equation*}
which due to Lemma \ref{lemma_Ll} coincides with \eqref{Ly0}.
$\Box$

Let us now prove the following useful expression for the partial derivative of $L_l$ with respect to $x$, while assuming that $\y$ does not depend on $x$.
\begin{lemma}
\label{lemma_Lx}
 Let $L_n$ be a polynomial defined by \eqref{Ll-notation} with $\Sigma_k$ being the sum from \eqref{Sigmak}. Let $\omega(Q_0)$ be defined by \eqref{omega_atQ0}. Assuming $x$ and $\y$ to be independent, we have
\begin{equation}
\label{Lx}
\frac{\partial L_n}{\partial x} = \frac{L_n}{2(x-\y)} - \frac{1}{2\omega(Q_0)}\frac{\partial^n}{\partial \y^n} \left\{ \frac{\omega(Q_0)}{x-\y} \right\}
=-\frac{n}{2\omega(Q_0)(x-\y)}\frac{\partial^{n-1}}{\partial \y^{n-1}} \left\{ \frac{\omega(Q_0)}{x-\y} \right\}\,.
\end{equation}
\end{lemma}
{\it Proof.} From definition \eqref{Ll-notation} of $L_n$ we have
\begin{multline*}
\frac{\partial L_n}{\partial x} =  \sum_{k=1}^n \sum_{\substack{p_1+2p_2+\dots +np_n=n \\ p_1,\dots, p_n\geq 0}}\frac{ n!\,\Sigma_1^{p_1} \Sigma_2^{p_2}\dots \Sigma_n^{p_n} p_k (-k)}{2^{\sum_{k=1}^n p_k} \prod_{k=1}^n p_k! \prod_{k=1}^n k^{p_k}\Sigma_k (x-\y)^{k+1}}
\\
=-\frac{1}{2}\sum_{k=1}^n \sum_{\substack{\hat p_1+2\hat p_2+\dots +n\hat p_n=n-k \\ \hat p_j=p_j, \mbox{ \small if } j\neq k, \,\,\, \hat p_k=p_k-1}}\frac{ n!\,\Sigma_1^{\hat p_1} \Sigma_2^{\hat p_2}\dots \Sigma_n^{\hat p_n} }{2^{\sum_{k=1}^n \hat p_k} \prod_{k=1}^n \hat p_k! \prod_{k=1}^n k^{\hat p_k}(x-\y)^{k+1}}
=-\frac{1}{2}\sum_{k=1}^n \frac{ n!\,L_{n-k} }{(n-k)!(x-\y)^{k+1}}\,.
\end{multline*}
On the other hand, by Leibniz rule and Lemma \ref{lemma_Ll},
\begin{multline}
\label{Leibniz}
\frac{\partial^n}{\partial \y^n} \left\{ \frac{\omega(Q_0)}{x-\y} \right\}
 =  \sum_{k=0}^n  \frac{ n! }{k!(n-k)!}  \frac{\partial^{n-k}}{\partial \y^{n-k}} \omega(Q_0) \frac{\partial^{k}}{\partial \y^{k}} \left\{ \frac{1}{x-\y} \right\}
=\sum_{k=0}^n  \frac{ n! }{(n-k)!}    \frac{ \omega(Q_0) L_{n-k}}{(x-\y)^{k+1}}
\\
= \frac{ \omega(Q_0) L_n}{x-\y}  - 2\omega(Q_0)\frac{\partial L_n}{\partial x}\,,
\end{multline}
where in the last equality we used the above result for $\frac{\partial L_n}{\partial x}\,.$ The last expression implies the first equality in \eqref{Lx}.

Alternatively, we can rewrite the expression obtained for $\partial L_n/\partial x$  by changing the summation index to $s:=k-1$
\begin{equation*}
\frac{\partial L_n}{\partial x} =  -\frac{n}{2(x-\y)}\sum_{s=0}^{n-1} \frac{ (n-1)!\,L_{n-1-s} }{(n-1-s)!(x-\y)^{s+1}}\,,
\end{equation*}
which can be interpreted using the $(n-1)$th derivative of $\omega(Q_0)/(x-\y)$ as in the second equality of \eqref{Lx}.
$\Box$

\section{Rauch variation}
\label{sect_Rauch}

Recall that our goal is to obtain a differential equation for the projection of $Q_0$ on the $u$-sphere denoted by $\y$ as a function of $x$. To this end, we will start by differentiating relations \eqref{periods} with respect to $x$. To differentiate with respect to the branch point $x$, we use the Rauch variational formulas,  giving derivatives of $W$ with respect to simple branch points, see see \cite{Rauch, Fay92, KokoKoro}. To introduce the Rauch formulas, we need to consider a family of elliptic curves written in a form slightly different from \eqref{curve}. Namely, consider a family of compact surfaces corresponding to the curves
$$\{(u,v)\in\mathbb C^2 | v^2=(u-u_1)(u-u_2)(u-u_3)\}$$
parameterized by distinct values of the branch points $u_1, \, u_2, \, u_3$ varying in some small neighbourhood so that a homology basis $\{a,b\}$ can be chosen for all surfaces in the family in a way that $u(a)$ and $u(b)$ be independent of the branch points. Let $W(P, Q)$ be the fundamental Riemann bidifferential on defined on a  surface from this family with respect to the chosen homology basis. Let $\omega(P)=\frac{1}{2\pi \i}\oint_bW(P,Q)$ be the holomorphic normalized differential and $\tau=\oint_b\omega$ be the modulus of the elliptic surface. Then using notation (\ref{evaluation}) for evaluation of differentials at a ramification point $P_{u_j}$ with respect to the standard local parameter $\zeta_{u_j}(P)=\sqrt{u(P)-u_j}$,  for $j=1,2,3$ we have (see Theorem 5 of \cite{KokoKoro} or (2.14) in \cite{KokoKoro2}):
\begin{equation}
\label{RauchWgeneral}
\frac{\partial^{\rm Rauch} }{\partial  u_j} W(P,Q)= \frac{1}{2} W(P,P_{u_j}) W(P_{u_j},Q)\,,
\end{equation}
where the notation $\partial^{\rm Rauch}$ means that  the projections of the points $P$ and $Q$ on the $u$-sphere are kept fixed for the differentiation with respect to $u_j$. By performing the integration over the $b$ cycles, we deduce from \eqref{RauchWgeneral} the variational formulas for $\omega$ and $\tau$, see (3.21) in \cite{Fay92} or (2.12), (2.13) in \cite{KokoKoro2}:
\begin{equation}
\label{Rauchomega}
\frac{\partial^{\rm Rauch}  \omega(P)}{ \partial u_j} = \frac{1}{2} W(P,P_{u_j}) \omega(P_{u_j});
\end{equation}
\begin{equation}
\label{Rauchtau}
\frac{\partial^{\rm Rauch}  {\tau}}{\partial  {u_j}} = \pi {\rm i}\, \omega^2(P_{u_j}),
\end{equation}
We will also need to differentiate with respect to $x$ various quantities evaluated at $P_x$. The above Rauch formulas cannot be used to achieve this, thus we need the following lemma.
\begin{lemma} Let $\surf$ be the compact Riemann surfaces of genus one corresponding to the algebraic curves \eqref{curve} with $x\in \X\subset \mathbb C\setminus\{0,1\}$, where the set $X$ is small enough so that the homology basis $\{a, b\}$ can be chosen on each surface $\surf$ in a way that the projections of the basis cycles on the $u$-sphere, $u(a)$ and $u(b)$ are independent of $x\in \X\,.$ Let the notation for local parameters and evaluation of differentials at a point be defined by \eqref{coordinates} and \eqref{evaluation}.  Let $\omega$ be the holomorphic normalized differential \eqref{omega}, \eqref{I0} and $W(P,Q)$ be the fundamental Riemann bidifferential \eqref{W} on $\surf\,.$ Then the following variational formulas hold:
\label{lemma_epsilon}
\begin{equation*}
\frac{d\omega(P_x)}{dx} = -\frac{1}{2} \omega(P_0)W(P_x, P_0)-\frac{1}{2} \omega(P_1)W(P_x, P_1);
\end{equation*}
\begin{equation*}
\frac{dW(Q_0,P_x)}{dx} = -\frac{1}{2}W(Q_0,P_0)W(P_0,P_x)-\frac{1}{2}W(Q_0,P_1)W(P_1,P_x) -\frac{\partial}{\partial\y} W( Q_0,P_x)+\frac{\partial}{\partial\y} W( Q_0,P_x)\y'\,.
\end{equation*}
\end{lemma}
{\it Proof.} Let $\varepsilon$ be a complex number sufficiently close to zero and let  the ramified covering  $u_\varepsilon: \surf \to \mathbb CP^1$ be defined by $u_\varepsilon: P \mapsto u(P)+\varepsilon\,,$ where $P=(u(P),v(P))$ is a point of the curve \eqref{curve}. For $\varepsilon=0$ we obtain the covering $u$ from Section \ref{sect_diff}.

The covering $u_\varepsilon$  is ramified at the  ramification points $P_0=(0,0)$, $P_1=(1,0)$, $P_x=(x,0)$ and $P_\infty=(\infty, \infty)$ for any $\varepsilon$;  the branch points of $u_\varepsilon$ are at $u_1:=\varepsilon, \,u_2:=1+\varepsilon,\, u_3:=x+\varepsilon,\, \infty$.

Denote by $W_\varepsilon(Q_0, P_x)$ the Riemann bidifferential $W$ \eqref{W} evaluated at the points $Q_0$ and $P_x$ according to definition \eqref{evaluation} with respect to the standard local parameters induced by the covering $u_\varepsilon\,.$ Assume for a moment that the point $Q_0$ is an arbitrary regular point of the covering, in particular that its position does not depend on $x$. Note that
 $W^\varepsilon(Q_0, P_x)$ is a function of the branch points: $u_1^\varepsilon:=\varepsilon, \,u_2^\varepsilon:=1+\varepsilon,\, u_3^\varepsilon:=x+\varepsilon$ and of $\y^\varepsilon:=\y+\varepsilon,$ the $u_\varepsilon$ coordinate of $Q_0$. Thus we have
 \begin{equation}
 \label{epsilon_temp}
 \frac{d}{d\varepsilon} W^\varepsilon(Q_0, P_x) = \frac{\partial W^\varepsilon(Q_0, P_x)}{\partial u_1^\varepsilon}+\frac{\partial W^\varepsilon(Q_0, P_x)}{\partial u_2^\varepsilon}+\frac{\partial W^\varepsilon(Q_0, P_x)}{\partial u_3^\varepsilon}+\frac{\partial W^\varepsilon(Q_0, P_x)}{\partial \y^\varepsilon}\,.
 \end{equation}
On the other hand, the standard local parameters \eqref{coordinates} at $Q_0$ and at $P_x$ induced by the covering $u_\varepsilon$ with respect to which the evaluation $W^\varepsilon(Q_0, P_x)$ is computed do not depend on $\varepsilon:$ we have $\zeta_{u_3}(P)=\sqrt{u_\varepsilon(P)-u_3}=\sqrt{u(P)-x}=\zeta_x(P)$ and $\zeta_{Q_0}(P)=u_\varepsilon(P)-u_\varepsilon(Q_0) = u(P)-u(Q_0)$ and thus
\begin{equation}
\label{Wpullback}
W^\varepsilon(Q_0,P_x)=\frac{W(P, Q)}{d\zeta_{u_3}(P)d\zeta_{Q_0}(Q) }\Big{|}_{\underset{Q=Q_0}{P=P_x}} = \frac{W(P, Q)}{d\zeta_{x}(P)d\zeta_{Q_0}(Q) }\Big{|}_{\underset{Q=Q_0}{P=P_x}} = W(Q_0,P_x)\,.
\end{equation}
This implies that the derivative with respect to $\varepsilon$ in \eqref{epsilon_temp} vanishes and, evaluating \eqref{epsilon_temp} at $\varepsilon=0\,,$ we have
 \begin{equation}
 \label{epsilon_temp1}
\frac{\partial W(Q_0, P_x)}{\partial u_1^0}+\frac{\partial W(Q_0, P_x)}{\partial u_2^0}+\frac{\partial W(Q_0, P_x)}{\partial u_3^0}+\frac{\partial W(Q_0, P_x)}{\partial \y}=0\,.
 \end{equation}
 Now, the derivatives with respect to $u_1^0$ and $u_2^0$ can be obtained by Rauch formulas \eqref{RauchWgeneral}. Using at the same time $u_0^\varepsilon=x,$ we rewrite \eqref{epsilon_temp1} in the form:
 \begin{equation}
 \label{epsilon_temp2}
\frac{\partial W(Q_0, P_x)}{\partial x}=-\frac{1}{2} W(Q_0, P_0)W(P_0, P_x)-\frac{1}{2}W(Q_0, P_1)W( P_1, P_x)-\frac{\partial W(Q_0, P_x)}{\partial \y}\,.
 \end{equation}
 Thus we have obtained the derivative of $W(Q_0, P_x)$ with respect to $x$ assuming that $Q_0$ does not depend on $x$. The full derivative with respect to $x$ can then be written as
 \begin{equation*}
\frac{d W(Q_0, P_x)}{d x}=\frac{\partial W(Q_0, P_x)}{\partial x} +\y'\frac{\partial W(Q_0, P_x)}{\partial \y}\,,
 \end{equation*}
 which, together with \eqref{epsilon_temp2}, gives the second formula of the lemma.
The derivative with respect to $x$ of $\omega(P_x)$ is obtained in a similar although simpler way, as there is no dependence on $Q_0\,.$
$\Box$

\begin{corollary}
\label{cor_epsilon}
Let the normalizing constant $I^x$ be defined by \eqref{WPPx} and let  $L_1$ stand, as before, for $L_1=-\frac{1}{2} \left(\frac{1}{\y} + \frac{1}{\y-1} + \frac{1}{\y-x}\right) $ according to \eqref{L}, \eqref{Sigmak} and \eqref{Ll-notation}.
Then, with the assumptions and notation of Lemma \ref{lemma_epsilon}, the variational formulas of the lemma can be rewritten as follows:
\begin{equation}
\label{domega}
\frac{d\omega(P_x)}{dx} = \frac{\omega(P_x)}{2}\left(I^x \omega(P_x)-   \frac{1}{x} - \frac{1}{x-1} \right),
\end{equation}
and
\begin{multline}
\label{dW}
\frac{dW(Q_0,P_x)}{dx} = -\frac{\omega(Q_0)}{2\omega(P_x)} \left(  \frac{1}{\y-x}\left( \frac{1}{\y} +\frac{1}{\y-1}- \frac{1}{x} - \frac{1}{x-1}\right) + \frac{1}{x} - \frac{1}{x-1} \right)+ \frac{1}{2} (I^x)^2\omega(Q_0) \omega(P_x)
\\
-\frac{\omega(Q_0)I^x}{2} \left(  \frac{1}{\y} +\frac{1}{\y-1} + \frac{1}{x} + \frac{1}{x-1}  \right)
%\\
+(\y'-1) \left(  -\frac{\omega(Q_0)}{\omega(P_x)(\y-x)^2} + W(P_x, Q_0)L_1  \right).
\end{multline}
\end{corollary}
{\it Proof.} The derivative of $\omega(P_x)$ is obtained by plugging expressions \eqref{W0x} and \eqref{W1x} for $W(P_0,P_x)$ and $W(P_1,P_x)$ in terms of $I^x$ into the expression of Lemma \ref{lemma_epsilon} and then using relation \eqref{zero} as well as explicit expressions \eqref{omega_atP0}-\eqref{omega_atPx} for $\omega(P_j)\,.$

For the derivative of $W(Q_0, P_x)$, we rewrite the expression from Lemma \ref{lemma_epsilon} plugging in expressions \eqref{W0x}-\eqref{W1x} and \eqref{WPP0}-\eqref{WPPx} with $P=Q_0\,:$
\begin{multline*}
\frac{dW(Q_0,P_x)}{dx} = - \frac{1}{2} \left( \frac{1}{\omega(P_0)\y}+I^0\right)\omega(Q_0)\left( I^x-\frac{1}{\omega(P_x)x}\right)\omega(P_0)
\\
-\frac{1}{2}\left( \frac{1}{\omega(P_1)(\y-1)}+I^1\right)\omega(Q_0)\left( I^x-\frac{1}{\omega(P_x)(x-1)}\right)\omega(P_1)
\\
+(\y'-1) \frac{\partial}{\partial\y}\left\{\left( \frac{1}{\omega(P_x)(\y-x)}+I^x\right)\omega(Q_0) \right\}.
\end{multline*}
Now, we express $I^0$ and $I^1$ from \eqref{I0x} and \eqref{I1x}, evaluating at the same time the derivative with respect to $\y$ explicitly differentiating expression \eqref{WPPx} with $P=Q_0$ for $W(Q_0, P_x)\,.$ This yields:
\begin{multline*}
\frac{dW(Q_0,P_x)}{dx} = - \frac{\omega(Q_0)}{2} \left( \frac{1}{\y}+\frac{\omega^2(P_0)I^x}{\omega(P_x)} - \frac{\omega^2(P_0) I_0^2 x}{4}\right)\left( I^x-\frac{1}{\omega(P_x)x}\right)
\\
-\frac{\omega(Q_0)}{2}\left( \frac{1}{\y-1}+\frac{\omega^2(P_1)I^x}{\omega(P_x)} - \frac{(x-1)I_0^2 \omega^2(P_1)}{4}\right)\left( I^x-\frac{1}{\omega(P_x)(x-1)}\right)
\\
+(\y'-1) \left(- \frac{\omega(Q_0)}{\omega(P_x)(\y-x)^2} + W(P_x, Q_0)L_1 \right).
\end{multline*}
Now we transform this expression multiplying terms in the first two lines and noting that due to \eqref{I0x} $\frac{\omega^2(P_0) I_0^2 x}{4}=1$ and due to \eqref{I1x} $\frac{(x-1)I_0^2 \omega^2(P_1)}{4}=-1$. Then using again \eqref{omega_atP0}-\eqref{omega_atPx} to obtain $\frac{\omega^2(P_0)}{x\omega^2(P_x)}=1-\frac{1}{x}$ and $\frac{\omega^2(P_1)}{\omega^2(P_x)(x-1)}=-1-\frac{1}{x-1}$ while collecting terms with equal powers of $I^x$ and using \eqref{zero}, we obtain the claim of the lemma.
$\Box$

\section{Isoperiodicity equation}
\label{sect_equation}
Let us first rewrite the expression for the $b$-period of $\Omega$ from \eqref{periods} more explicitly using the Bell polynomial $L_n.$ Given that $\omega^{(n)}(Q_0)$ is the $n$th term in the Taylor expansion \eqref{omega-Taylor} it is proportional to the $n$-fold derivative of the function $\omega(P)/d\zeta_{Q_0}(P)$ evaluated at $P=Q_0,$ which is the same as the $n$-fold derivative with respect to $\y$ of the function $\omega(Q_0)$ given by \eqref{omega_atQ0}. Due to Lemma \ref{lemma_Ll} and notation \eqref{Ll-notation}, we have
\begin{equation}
\label{periodL}
\bb=\oint_b\Omega=\frac{2\pi\i}{n+1} \omega^{(n)}(Q_0)+\aa\tau=\frac{2\pi\i}{(n+1)!} \omega(Q_0)L_n+\aa\tau\,.
\end{equation}

Differentiating this equality with respect to $x$, we can express $\y'=d\y/dx$ as follows.
\begin{lemma}
\label{lemma_y0'}
Let $\Omega$ be the meromorphic differential of the second kind defined by \eqref{Omega} on the each compact surface $\surf$ corresponding to the curve \eqref{curve} for $x\in \X\subset \mathbb C\setminus\{0,1\}$. Let the position $Q_0\in\surf$ of the pole of $\Omega$ be defined by the conditions $\oint_a\Omega= {\rm const}$ and  $\oint_b\Omega= {\rm const}.$ Then the $u$-coordinate $\y=u(Q_0)$ of the pole as a function of $x$ satisfies
\begin{equation}
\label{y0'}
\y'=-\frac{1}{2L_{n+1}\omega(Q_0)}\left( \omega(P_x)W(Q_0,P_x) L_n -  \frac{n}{(x-\y)}\frac{\partial^{n-1}}{\partial \y^{n-1}} \left\{ \frac{\omega(Q_0)}{x-\y} \right\}    +\aa(n+1)! \omega^2(P_x)\right)\,.
\end{equation}
\end{lemma}
{\it Proof.} Given that the $b$-period of $\Omega$ \eqref{Omega} is assumed constant, the derivative with respect to $x$ of the right-hand side in \eqref{periodL} vanishes. Thus differentiation of \eqref{periodL}  with respect to the branch point $x$ leads to
\begin{multline*}
0=\frac{2\pi\i}{(n+1)!} \frac{d\omega(Q_0)}{dx}L_n+\frac{2\pi\i}{(n+1)!} \omega(Q_0)\frac{dL_n}{dx}+\aa\frac{d\tau}{dx}
\\
=\frac{2\pi\i}{(n+1)!} \left( \frac{\partial^{\rm Rauch}  \omega(Q_0)}{ \partial x} +\frac{\partial \omega(Q_0)}{\partial \y} \y' \right)L_n +\frac{2\pi\i}{(n+1)!} \omega(Q_0)\left(\frac{ \partial L_n}{\partial x} +\frac{ \partial L_n}{\partial \y}\y' \right) +\aa\frac{\partial^{\rm Rauch}\tau}{\partial x}\,.
\end{multline*}
Plugging in Rauch derivatives \eqref{Rauchomega}, \eqref{Rauchtau} and \eqref{Ll} from Lemma \ref{lemma_Ll} as well as \eqref{Ly0} and the second equality of Lemma \ref{lemma_Lx} for derivatives of $L_n,$ we obtain the claim of the lemma.
$\Box$

\begin{example}
\label{example_y0'-n=0}
In the case $n=0$, the differential $\Omega$ \eqref{Omega} is given by
\begin{equation*}
\Omega(P)=W(P, Q_0)+\aa \omega(P)\,.
\end{equation*}
Applying Lemma \ref{lemma_y0'}, and using $L_0=1$ and $L_1(-\Sigma_1)=\frac{1}{2}\Sigma_1$ with $\Sigma_1= -\frac{1}{\y}-\frac{1}{\y-1}-\frac{1}{\y-x}\,,$ we can express $\y'$ in terms of $\Omega(P_x)$ as follows:
\begin{equation*}
\y'= \frac{ \omega(P_x)\Omega(P_x)}{\omega(Q_0)\left(\frac{1}{\y}+\frac{1}{\y-1}+\frac{1}{\y-x} \right)}\,.
\end{equation*}
\end{example}

Now, our aim is to differentiate the expression for $\y'$ from Lemma \ref{lemma_y0'} the second time with respect to the branch point $x$ in order to  derive a second order ODE for the function $\y(x)\,.$
\begin{theorem}
\label{thm_ode}
Let $n\geq 0$ be a fixed integer and let $\Omega$ be the meromorphic differential of the second kind defined by \eqref{Omega} on each compact surface $\surf$ corresponding to the curve \eqref{curve} for $x\in \X\subset \mathbb C\setminus\{0,1\}$. The set $\X$ is assumed simply connected to allow for a choice of homology basis $\{a,b\}$ such that $u(a)$ and $u(b)$ be fixed for all $x\in \X$. Let the position $Q_0\in\surf$ of the pole of $\Omega$ be defined by the conditions $\oint_a\Omega= {\rm const}$ and  $\oint_b\Omega= {\rm const}.$ Then the $u$-coordinate $\y=u(Q_0)$ of the pole as a function of $x$ satisfies the following second order ordinary differential equation:
\begin{multline}
\label{ode}
\y''=-(\y')^2 \frac{L_{n+2}}{L_{n+1}}- \y'\left( \frac{1}{x}+\frac{1}{x-1}+ \frac{ 1}{\y-x}\right)
+ \y'   \frac{(n+1)! }{(x-\y)^{n+2}L_{n+1}} \sum_{s=0}^{n} \frac{L_s(x-\y)^s}{s!}
\\
+\frac{n! \left( \frac{1}{x} + \frac{1}{x-1} + \frac{2}{\y-x}\right)}{2(x-\y)^{n+1}L_{n+1}} \sum_{s=0}^{n-1}\frac{L_s(x-\y)^s}{s!}
-\frac{n!}{2(x-\y)^{n+2}L_{n+1}}\sum_{s=0}^{n-1}\frac{(n-s)L_s(x-\y)^s}{s!}
\\
-\frac{n!}{4(x-\y)^{n+2}L_{n+1}}\sum_{s=1}^{n-1}\sum_{k=0}^{s-1}\frac{L_k(x-\y)^k}{k!}
-\frac{L_n}{4L_{n+1}} \left( \frac{2}{\y-x} \left(  \frac{1}{x} + \frac{1}{x-1} \right) + \frac{3}{(\y-x)^2}+  \frac{1}{x-1}-\frac{1}{x} \right) \,.
\end{multline}
\end{theorem}

{\it Proof.} Let us differentiate expression \eqref{y0'} for $\y'$  with respect to $x$ in the following way:
\begin{multline}
\label{y''1}
\y''=-\frac{\frac{d}{dx}\left\{L_{n+1}\omega(Q_0))\right\}}{L_{n+1}\omega(Q_0)}\y'
\\
-\frac{1}{2L_{n+1}\omega(Q_0)}\frac{d}{dx}\left( \omega(P_x)W(Q_0,P_x) L_n -  \frac{n}{(x-\y)}\frac{\partial^{n-1}}{\partial \y^{n-1}} \left\{ \frac{\omega(Q_0)}{x-\y} \right\}    +\aa(n+1)! \omega^2(P_x)\right)
\,.
\end{multline}
Using \eqref{Ly0} and \eqref{Lx} to differentiate polynomial $L_{n+1}$ and \eqref{Ll} to differentiate $\omega(Q_0)$ with respect to $\y$ and the Rauch formula \eqref{Rauchomega} with $u_j=x$ to differentiate $\omega(Q_0)$ with respect to $x\,,$ the first line in \eqref{y''1} gives
\begin{multline*}
-\frac{\frac{d}{dx}\left\{L_{n+1}\omega(Q_0))\right\}}{L_{n+1}\omega(Q_0)}\y'=-\frac{\y'}{L_{n+1}}\left( \frac{\partial L_{n+1}}{\partial x} +\y'\frac{\partial L_{n+1}}{\partial \y}   \right)
-\frac{\y'}{\omega(Q_0)}\left(\frac{\partial^{\rm Rauch}}{\partial x}\Big{|}_{\y=const}\omega(Q_0)+\y'\frac{\partial \omega(Q_0)}{\partial \y}\right)
\\
=-\frac{\y'}{L_{n+1}}\left( -\frac{n+1}{2\omega(Q_0)(x-\y)}\frac{\partial^{n}}{\partial \y^{n}} \left\{ \frac{\omega(Q_0)}{x-\y} \right\} +\y'(L_{n+2}-L_1L_{n+1})   \right)
\\
-\frac{\y'}{\omega(Q_0)}\left(\frac{1}{2}\omega(P_x)W(P_x,Q_0)+\y' \omega(Q_0)L_1\right)
\end{multline*}

Differentiating the first term in the second line of \eqref{y''1} is done with the help of  formulas \eqref{domega} and \eqref{dW} from corollary \ref{cor_epsilon} for the full derivatives of $\omega(P_x)$ and of $W(Q_0, P_x)$ as well as formulas \eqref{Lx} and \eqref{Ly0} for derivatives of $L_n:$
\begingroup
\allowdisplaybreaks
\begin{multline*}
\frac{d}{dx}\left( \omega(P_x)W(Q_0,P_x) L_n\right)=\frac{1}{2}\left({I^x \omega(P_x)}-  \frac{1}{x} - \frac{1}{x-1} \right)\omega(P_x)W(Q_0,P_x) L_n
\\
+\left[  -\frac{\omega(Q_0)}{2\omega(P_x)} \left(  \frac{1}{\y-x}\left( \frac{1}{\y} - \frac{1}{x}\right) + \frac{1}{x} - \frac{1}{x-1}+ \frac{1}{\y-x}\left( \frac{1}{\y-1} - \frac{1}{x-1}\right) \right) \right.
\\
-\frac{\omega(Q_0)I^x}{2} \left(  \frac{1}{\y} +\frac{1}{\y-1} + \frac{1}{x} + \frac{1}{x-1}  \right) + \frac{1}{2} (I^x)^2\omega(Q_0) \omega(P_x)
\\
\left.
+(\y'-1) \left(  -\frac{\omega(Q_0)}{\omega(P_x)(\y-x)^2} + W(P_x, Q_0)L_1  \right)\right]\omega(P_x) L_n
\\
+\left( -\frac{n}{2\omega(Q_0)(x-\y)}\frac{\partial^{n-1}}{\partial \y^{n-1}} \left\{ \frac{\omega(Q_0)}{x-\y} \right\}+\y'(L_{n+1}-L_nL_1)\right)\omega(P_x)W(Q_0,P_x) \,.
\end{multline*}
\endgroup
Let us differentiate the second term in the second line of  \eqref{y''1}:
\begin{multline*}
\frac{d}{dx}\left( -  \frac{n}{(x-\y)}\frac{\partial^{n-1}}{\partial \y^{n-1}} \left\{ \frac{\omega(Q_0)}{x-\y} \right\}\right)
\\
=  \frac{n(1-\y')}{(x-\y)^2}\frac{\partial^{n-1}}{\partial \y^{n-1}} \left\{ \frac{\omega(Q_0)}{x-\y} \right\}
  -  \frac{n}{(x-\y)}\frac{\partial}{\partial x}\frac{\partial^{n-1}}{\partial \y^{n-1}} \left\{ \frac{\omega(Q_0)}{x-\y} \right\}
%  \\
  -  \frac{n\y'}{(x-\y)}\frac{\partial^{n}}{\partial \y^{n}} \left\{ \frac{\omega(Q_0)}{x-\y} \right\}\,.
\end{multline*}
To compute the derivative with respect to $x$ in the middle term, let us rewrite in terms of polynomials $L_n$ it using the Leibniz rule as in \eqref{Leibniz}:
\begin{multline*}
  \frac{\partial}{\partial x}\frac{\partial^{n-1}}{\partial \y^{n-1}} \left\{ \frac{\omega(Q_0)}{x-\y} \right\}
  =\frac{\partial}{\partial x} \sum_{s=0}^{n-1} \frac{(n-1)! L_s \omega(Q_0)}{s!(x-\y)^{n-s}}
\\
  =\frac{1}{2} \frac{\omega(P_x)W( P_x, Q_0)}{\omega(Q_0)}\frac{\partial^{n-1}}{\partial \y^{n-1}} \left\{ \frac{\omega(Q_0)}{x-\y} \right\}
%  \\
-  \sum_{s=0}^{n-1} \frac{(n-1)!(n-s) L_s \omega(Q_0)}{s!(x-\y)^{n-s+1}}
-\sum_{s=0}^{n-1} \frac{(n-1)! s\frac{\partial^{s-1}}{\partial \y^{s-1}} \left\{ \frac{\omega(Q_0)}{x-\y} \right\} }{2s!(x-\y)^{n-s+1}}\,,
\end{multline*}
where in the first term we differentiated $\omega(Q_0)$ according to Rauch \eqref{Rauchomega} and in the last term we used \eqref{Lx} to differentiate $L_s\,.$

The derivative of the last term in \eqref{y''1} is a straightforward application of \eqref{domega} from Corollary \ref{cor_epsilon}.

Collecting all these derivatives in \eqref{y''1} and cancelling two terms of the form $(\y')^2L_1$ and two terms of the form $\y'W(P_x, Q_0)\omega(P_x) L_nL_1$, we obtain
\begingroup
\allowdisplaybreaks
\begin{multline}
\label{y''2}
\y''=-\frac{\y'}{L_{n+1}}\left( -\frac{n+1}{2\omega(Q_0)(x-\y)}\frac{\partial^{n}}{\partial \y^{n}} \left\{ \frac{\omega(Q_0)}{x-\y} \right\} +\y'L_{n+2}   \right)
%\\
-\frac{\y'}{2}\frac{\omega(P_x)W(P_x,Q_0)}{\omega(Q_0)}
\\
-\frac{1}{2L_{n+1}\omega(Q_0)}\left\{ \frac{1}{2}\left({I^x \omega(P_x)}-  \frac{1}{x} - \frac{1}{x-1} \right)\omega(P_x)W(Q_0,P_x) L_n \right.
\\
+\omega(P_x) L_n\left[  -\frac{\omega(Q_0)}{2\omega(P_x)} \left(  \frac{1}{\y-x}\left( \frac{1}{\y} + \frac{1}{\y-1} - \frac{1}{x}- \frac{1}{x-1}\right) + \frac{1}{x} - \frac{1}{x-1} \right) \right.
\\
\left.
-\frac{\omega(Q_0)I^x}{2} \left(  \frac{1}{\y} +\frac{1}{\y-1} + \frac{1}{x} + \frac{1}{x-1}  \right) + \frac{1}{2} (I^x)^2\omega(Q_0) \omega(P_x)
   -\frac{\omega(Q_0)(\y'-1)}{\omega(P_x)(\y-x)^2} -    W(P_x, Q_0)L_1 \right]
\\
\left.
+\left( -\frac{n}{2\omega(Q_0)(x-\y)}\frac{\partial^{n-1}}{\partial \y^{n-1}} \left\{ \frac{\omega(Q_0)}{x-\y} \right\}+\y'L_{n+1}\right)\omega(P_x)W(Q_0,P_x) \right\}
\\
-\frac{1}{2L_{n+1}\omega(Q_0)}\left[    \frac{n(1-\y')}{(x-\y)^2}\frac{\partial^{n-1}}{\partial \y^{n-1}} \left\{ \frac{\omega(Q_0)}{x-\y} \right\} \right.
\\
  -  \frac{n}{(x-\y)} \left( \frac{1}{2} \frac{\omega(P_x)W( P_x, Q_0)}{\omega(Q_0)}\frac{\partial^{n-1}}{\partial \y^{n-1}} \left\{ \frac{\omega(Q_0)}{x-\y} \right\}
%  \\
-  \sum_{s=0}^{n-1} \frac{(n-1)!(n-s) L_s \omega(Q_0)}{s!(x-\y)^{n-s+1}}
-\sum_{s=0}^{n-1} \frac{(n-1)! s\frac{\partial^{s-1}}{\partial \y^{s-1}} \left\{ \frac{\omega(Q_0)}{x-\y} \right\} }{2s!(x-\y)^{n-s+1}}\right)
\\
\left.  -  \frac{n\y'}{(x-\y)}\frac{\partial^{n}}{\partial \y^{n}} \left\{ \frac{\omega(Q_0)}{x-\y} \right\}   \right]
%\\
-\frac{\aa(n+1)! \omega^2(P_x)}{2L_{n+1}\omega(Q_0)} \left( I^x \omega(P_x)-   \frac{1}{x} - \frac{1}{x-1} \right)
\,.
\end{multline}
\endgroup
Let us now note that the first terms in lines 5 and 7 coincide and their sum can be rewritten, expressing $W(Q_0, P_x)$ as in \eqref{WQ0x}, in the following form:
\begin{multline}
\label{tempy''}
\frac{n}{2L_{n+1}\omega(Q_0)(x-\y)}    \frac{\omega(P_x)W( P_x, Q_0)}{\omega(Q_0)}\frac{\partial^{n-1}}{\partial \y^{n-1}} \left\{ \frac{\omega(Q_0)}{x-\y} \right\}
\\
=
\frac{n}{2L_{n+1}\omega(Q_0)(x-\y)}   \left(  I^x\omega(P_x) -   \frac{1}{x} - \frac{1}{x-1} + \frac{1}{\y-x} +   \frac{1}{x} + \frac{1}{x-1} \right) \frac{\partial^{n-1}}{\partial \y^{n-1}} \left\{ \frac{\omega(Q_0)}{x-\y} \right\},
\end{multline}
where some terms were added and subtracted for the next step. Let us now combine the first three terms in \eqref{tempy''} and the last group of terms in \eqref{y''2}:
\begin{multline*}
 -\frac{1}{2L_{n+1}\omega(Q_0)} \left(  I^x\omega(P_x) -   \frac{1}{x} - \frac{1}{x-1}  \right) \left(
-\frac{n}{(x-\y)}   \frac{\partial^{n-1}}{\partial \y^{n-1}} \left\{ \frac{\omega(Q_0)}{x-\y} \right\} +\aa(n+1)! \omega^2(P_x) \right)
\\
=\left(  I^x\omega(P_x) -   \frac{1}{x} - \frac{1}{x-1}  \right)\left(\y'+\frac{\omega(P_x)W( P_x, Q_0)L_n}{2L_{n+1}\omega(Q_0)}\right),
\end{multline*}
where we used expression \eqref{y0'} for $\y'$ from Lemma \ref{lemma_y0'}. The last obtained expression cancels partially against the term in the second line of \eqref{y''2}. Incorporating all these transformations in \eqref{y''2}, opening some of the parentheses,
replacing $\frac{\omega(P_x)W(P_x,Q_0)}{\omega(Q_0)}$ by $\frac{1}{\y-x} + I^x\omega(P_x)$ according to \eqref{WQ0x}, we obtain:
\begin{multline}
\label{y''3}
\y''=-\frac{\y'}{L_{n+1}}\left( -\frac{n+1}{2\omega(Q_0)(x-\y)}\frac{\partial^{n}}{\partial \y^{n}} \left\{ \frac{\omega(Q_0)}{x-\y} \right\} +\y'L_{n+2}   \right)
-\y'\left( \frac{1}{\y-x} + I^x\omega(P_x)\right)
\\
+
 \frac{ L_n}{4L_{n+1}} \left(  \frac{1}{\y-x}\left( \frac{1}{\y} + \frac{1}{\y-1} - \frac{1}{x}- \frac{1}{x-1}\right) + \frac{1}{x} - \frac{1}{x-1} \right)
\\
+\frac{I^x\omega(P_x) L_n}{4L_{n+1}} \left(  \frac{1}{\y} +\frac{1}{\y-1} + \frac{1}{x} + \frac{1}{x-1}  \right)
- \frac{\omega^2(P_x)  L_n(I^x)^2}{4L_{n+1}}
   +\frac{(\y'-1) L_n}{2(\y-x)^2 L_{n+1}}
\\
+ \frac{ L_nL_1 }{2L_{n+1}}\left( \frac{1}{\y-x} + I^x\omega(P_x)\right)
   -\frac{n(1-\y')}{2L_{n+1}\omega(Q_0)(x-\y)^2}\frac{\partial^{n-1}}{\partial \y^{n-1}} \left\{ \frac{\omega(Q_0)}{x-\y} \right\}
\\
-  \frac{1}{2L_{n+1}}\sum_{s=0}^{n-1} \frac{n!(n-s) L_s }{s!(x-\y)^{n-s+2}}
-\frac{1}{4L_{n+1}\omega(Q_0)}\sum_{s=0}^{n-1} \frac{n! s\frac{\partial^{s-1}}{\partial \y^{s-1}} \left\{ \frac{\omega(Q_0)}{x-\y} \right\} }{s!(x-\y)^{n-s+2}}
\\
  +  \frac{1}{2L_{n+1}\omega(Q_0)}\frac{n\y'}{(x-\y)}\frac{\partial^{n}}{\partial \y^{n}} \left\{ \frac{\omega(Q_0)}{x-\y} \right\}
+\left(  I^x\omega(P_x) -   \frac{1}{x} - \frac{1}{x-1}  \right)\left(\y'+\frac{L_n}{4L_{n+1}}\left( \frac{1}{\y-x} + I^x\omega(P_x)\right)\right)
\\
+\frac{n}{2L_{n+1}\omega(Q_0)(x-\y)}   \left(   \frac{1}{\y-x} +   \frac{1}{x} + \frac{1}{x-1} \right) \frac{\partial^{n-1}}{\partial \y^{n-1}} \left\{ \frac{\omega(Q_0)}{x-\y} \right\}
\,.
\end{multline}
It remains to collect the terms with $I^x$ and $(I^x)^2$ to see that they vanish and note that the remaining terms are rational. We have

\begin{multline}
\label{y''6}
\y''=-\frac{\y'}{L_{n+1}}\left( -\frac{n+1}{2\omega(Q_0)(x-\y)}\frac{\partial^{n}}{\partial \y^{n}} \left\{ \frac{\omega(Q_0)}{x-\y} \right\} +\y'L_{n+2}   \right)
- \frac{\y'}{\y-x}
\\
+
 \frac{ L_n}{4L_{n+1}} \left(  \frac{1}{\y-x}\left( \frac{1}{\y} + \frac{1}{\y-1} - \frac{1}{x}- \frac{1}{x-1}\right) + \frac{1}{x} - \frac{1}{x-1} \right)
\\
   +\frac{(\y'-1) L_n}{2(\y-x)^2 L_{n+1}}
%\\
- \frac{ L_nL_1 }{2L_{n+1}(x-\y)}
   -\frac{n(1-\y')}{2L_{n+1}\omega(Q_0)(x-\y)^2}\frac{\partial^{n-1}}{\partial \y^{n-1}} \left\{ \frac{\omega(Q_0)}{x-\y} \right\}
\\
-  \frac{1}{2L_{n+1}}\sum_{s=0}^{n-1} \frac{n!(n-s) L_s }{s!(x-\y)^{n-s+2}}
-\frac{1}{4L_{n+1}\omega(Q_0)}\sum_{s=0}^{n-1} \frac{n! s\frac{\partial^{s-1}}{\partial \y^{s-1}} \left\{ \frac{\omega(Q_0)}{x-\y} \right\} }{s!(x-\y)^{n-s+2}}
\\
  +  \frac{1}{2L_{n+1}\omega(Q_0)}\frac{n\y'}{(x-\y)}\frac{\partial^{n}}{\partial \y^{n}} \left\{ \frac{\omega(Q_0)}{x-\y} \right\}
%\\
-\left(    \frac{1}{x} + \frac{1}{x-1}  \right)\y'
+\frac{L_n}{4L_{n+1}(x-\y)}\left(     \frac{1}{x} + \frac{1}{x-1}  \right)
\\
+\frac{n}{2L_{n+1}\omega(Q_0)(x-\y)}   \left(   \frac{1}{\y-x} +   \frac{1}{x} + \frac{1}{x-1} \right) \frac{\partial^{n-1}}{\partial \y^{n-1}} \left\{ \frac{\omega(Q_0)}{x-\y} \right\}
\,.
\end{multline}
Rewriting derivatives of $\frac{\omega(Q_0)}{x-\y}$ with respect to $\y$ in terms of polynomials $L_l$ as in \eqref{Leibniz}, we see that the obtained
equation is equivalent to the claim of the theorem.
$\Box$

\begin{example}
\label{example_n0}
For $n=0,$ the differential equation from Theorem \ref{thm_ode} reads
\begin{multline*}
\y'' =  \frac{(\y')^2}{2}  \left(\frac{1}{\y}+\frac{1}{\y-1}+\frac{1}{\y-x}\right)
 +\frac{(\y')^2\left( \frac{1}{\y^2}+\frac{1}{(\y-1)^2}+\frac{1}{(\y-x)^2}\right)}{\frac{1}{\y}+\frac{1}{\y-1}+\frac{1}{\y-x}}
\\
-\y'\left(  \frac{1}{x}+\frac{1}{x-1}+ \frac{ 1}{\y-x}  \right)
-   \frac{ 2\y'}{\left(\frac{1}{\y}+\frac{1}{\y-1}+\frac{1}{\y-x}\right)(\y-x)^2}
\\
+ \frac{  \frac{2}{x(\y-x)}+\frac{2}{(x-1)(\y-x)}+ \frac{1}{x-1}- \frac{1}{x} + \frac{3}{(\y-x)^2}}{2\left(\frac{1}{\y}+\frac{1}{\y-1}+\frac{1}{\y-x}\right)}
 \,.
\end{multline*}
\end{example}
\begin{example}
\label{example_n1}
For $n=1,$ the differential equation from Theorem \ref{thm_ode} reads
\begin{multline*}
\y'' =  \frac{(\y')^2}{2}  \left(\frac{1}{\y}+\frac{1}{\y-1}+\frac{1}{\y-x}\right)-\y'\left(  \frac{1}{x}+\frac{1}{x-1}+ \frac{ 1}{\y-x}  \right)
\\+\frac{1}{{\left(\frac{1}{\y}+\frac{1}{\y-1}+\frac{1}{\y-x}\right)^2+2\left(\frac{1}{\y^2}+\frac{1}{(\y-1)^2}+\frac{1}{(\y-x)^2}\right)}} \times
\\
\times \left[ 2(\y')^2\left(\left(\frac{1}{\y}+\frac{1}{\y-1}+\frac{1}{\y-x}\right)\left(\frac{1}{\y^2}+\frac{1}{(\y-1)^2}+\frac{1}{(\y-x)^2}\right)+2\left(\frac{1}{\y^3}+\frac{1}{(\y-1)^3}+\frac{1}{(\y-1)^3}\right)\right)  \right.
\\
+   \frac{4\y'}{(\y-x)^3}\left(\frac{x}{\y}+\frac{x-1}{\y-1}-1\right)
+\frac{2}{(x-\y)^2}\left(\frac{1}{x}+\frac{1}{x-1}+\frac{3}{\y-x}\right)
\\
\left.
+\frac{1}{2}\left(\frac{1}{\y}+\frac{1}{\y-1}+\frac{1}{\y-x}\right) \left( \frac{2}{\y-x} \left(  \frac{1}{x} + \frac{1}{x-1} \right) + \frac{3}{(\y-x)^2}+  \frac{1}{x-1}-\frac{1}{x} \right)\right] \,.
\end{multline*}
\end{example}
Theorem \ref{thm_ode} gives a necessary condition for the deformations of an elliptic curve carrying a differential $\Omega$ \eqref{Omega} to be isoperiodic in the sense of Definition \ref{def-iso}. The next theorem gives a sufficient condition.
\begin{theorem}
\label{thm_unicity}
Let $n\geq 0$ be a fixed integer, $x_0\in\mathbb C\setminus\{0,1\}$, and let $\Omega(x_0)$ be the meromorphic differential of the second kind defined by \eqref{Omega} on a compact Riemann surface $\mathcal L(x_0)$ of the curve \eqref{curve} with $x=x_0$ and some choice of a point $Q_0(x_0)\in\mathcal L(x_0)$ different from ramification points. Assume that $x_0$ and $Q_0(x_0)$ are such that $\omega^{(n+1)}(Q_0(x_0))\neq 0$ for $\omega^{(n+1)}$ defined by \eqref{omega-Taylor}. Let $\{a,b\}$ be a canonical homology basis on $\mathcal L(x_0)$ such that the projections $u(a)$ and $u(b)$ do not intersect a certain neighbourhood $\hat\X$ of $x_0$. Then there exists a unique continuous family  $(\surf, \Omega(x))$, $x\in\X\subset\hat\X,$ where the family of curves $\surf$  \eqref{curve}  admits a continuous section $Q_0(x)\in\surf$ such that
 $\Omega(x)$ is a differential \eqref{Omega} of the second kind with a pole of order $n+1$ at $Q_0(x)$ having $a$- and $b$-periods independent of $x\in\X$.
\end{theorem}
{\it Proof.}  Given that the $a$-period of $\Omega=\Omega(x)$ is constant by definition \eqref{Omega}, an isoperiodic deformation is defined by the condition $\oint_b\Omega(x)={\rm const}. $ The quantity $\oint_b\Omega(x)$ is a function of $x$ and $y_0.$  By the implicit function theorem, the relation $\oint_b\Omega(x)={\rm const}$ defines a function $y_0(x)$ in a neighbourhood of $x=x_0$ if $\frac{\partial}{\partial y_0}|_{x=x_0}\oint_b\Omega(x)\neq 0.$
Note first that due to \eqref{Omega} and \eqref{pole},  $\frac{\partial \Omega}{\partial y_0} = \frac{\partial W^{(n)}(P, Q_0)}{\partial y_0}$ is a differential of the second kind with a pole at $Q_0$ of the order one higher than the order of the pole of $W^{(n)}(P, Q_0)$ with the Laurent series at the pole of the form
\begin{equation*}
W^{(n)}(P, Q_0)\underset{P\sim Q_0}{=}\left( \frac{n+2}{(u(P)-\y)^{n+3}} + \mathcal O (u(P)-\y) \right)du(P)\,.
\end{equation*}
Moreover, the $a$-period of $\frac{\partial \Omega}{\partial y_0}$ vanishes just like the $a$-period of $\Omega$. This, due to the unicity of a normalized differential of the second kind with a given behaviour at the pole, implies
\begin{equation*}
\frac{\partial \Omega}{\partial y_0}=(n+2)W^{(n+1)}(P, Q_0).
\end{equation*}

Now, by \eqref{pole} and \eqref{bn} for $x\in\hat\X$ where the contour $b$ is fixed by the condition $u(b)=const,$ we have
\begin{equation*}
\frac{\partial}{\partial y_0}\oint_b\Omega = \oint_b\frac{\partial \Omega}{\partial y_0}=(n+2)\oint_b W^{(n+1)}(P, Q_0) = 2\pi\i \, \omega^{(n+1)}(Q_0(x_0)).
\end{equation*}
This is non-zero given the choice of $x_0$.
This proves the existence of an isoperiodic deformation in some neighbourhood $\X\subset\hat \X$ of $x_0.$

Now, Theorem \ref{thm_ode} states that for every continuous isoperiodic deformation, $y_0(x)$ should satisfy equation \eqref{ode}. In addition, we know that for an isoperiodic deformation, $y_0'(x)$ is given by \eqref{y0'}. There is a unique solution to \eqref{ode}, \eqref{y0'} passing through the point $y_0(x_0)$. This proves the unicity of an isoperiodic deformation of the pair $(\mathcal L, \Omega(x_0)).$
$\Box$
\\

Note that the condition given in Theorem \ref{thm_unicity} for $x_0$ is valid for a generic point $x_0\in\mathbb C\setminus\{0,1\}$ as $\omega^{(n+1)}(Q_0(x_0))=0$ for a set of measure zero in the space of elliptic curves.

\section{Applications: isoperiodic deformations of solutions to the Boussinesq equation}
\label{sect_application}

In this section, we assume that $A=0$ in \eqref{Omega}. We consider a genus one Riemann surface $\surf$ corresponding to the elliptic curve of equation \eqref{curve} with an arbitrary point $Q_0$ and the normalized differentials of the second kind
\begin{equation}\label{eq:normsecond}
\Omega^{(n+1)}_{Q_0}(P):=W^{(n)}(P, Q_0)\,,\quad P\in \surf,
\end{equation}
with a unique pole at $Q_0$ of order $n+2$.
Denote further
\begin{equation}\label{eq:UV}
U:=-\frac{1}{2\pi\i}\oint_b\Omega^{(1)}_{Q_0},\qquad V:=-\frac{1}{\pi\i}\oint_b\Omega^{(2)}_{Q_0}.
\end{equation}

 \begin{remark}
\label{rem:realell}
 If an elliptic curve  is given by \eqref{curve} with real $x$, then it admits an antiholomorphic involution
$$
\rho(u, v)=(\bar u, \bar v),
$$
see e.g. \cite{DubNat}; the canonical basis of cycles can then be chosen so that  $\rho a=a$, $\rho b=-b$. In this case, for the holomorphic normalized differential \eqref{omega}, the normalization coefficient $I_0$
is real and the $b$-period $\tau$ \eqref{eq:tau} is purely imaginary,  see e.g. \cite{DubNat}.
Thus, if $y_0=u(Q_0)$ is real and such that $y_0(y_0-1)(y_0-x)>0$, then both $U$ and $V$ \eqref{eq:UV} are real, see \eqref{eq:normsecond} and \eqref{bn}.
\end{remark}

We follow \cite{Dub1981} closely in this section while keeping our normalization \eqref{omega}, \eqref{I0}.
Consider a meromorphic function $\mu$ with a unique pole of order three at $Q_0$. Such a function exists for any point $Q_0$ of an elliptic
curve, as one can take an appropriate shift of the derivative $\wp'$ of the corresponding Weierstrass function. Due to the existence of the function $\mu$, for a surface of genus one, solutions of the Kadomtsev-Petviashvili equation (KP) in terms of the theta-function obtained from the Baker-Akhiezer function, see for example \cite{Dub1981},  reduce to solutions of the Boussinesq equation.
The Baker-Akhiezer function in question, denoted by $\psi$, is constructed from the polynomial
$$
q(k)=kX+k^2Y,
$$
where $k$ is the reciprocal of the local parameter $\zeta$ around $Q_0$ on $\surf$. We may assume that $\mu=k^3+\mathcal O(1)$ for $P$ close to $Q_0$.
The Baker-Akhiezer function has the following expansion at $Q_0$:
$$
\psi(X, Y;P)=e^{kX+k^2Y}\big(1+\frac{\xi_1}{k}+\frac{\xi_2}{k^2}+\dots \big), \qquad P\sim Q_0,
$$
where $\xi_i$, $i=1,2$ are functions of $X$ and $Y$ to be determined below.
Consider the operators
$$
\mathbb L=\partial^2_{X}+\u, \quad \mathbb A=\partial^3_{X}+\frac{3}{2}\partial_X + \w,
$$
where
$$
\u=-2\frac{\partial \xi_1}{\partial X},\quad \w=3\xi_1\frac{\partial \xi_1}{\partial X}+3\frac{\partial^2 \xi_1}{\partial X^2}-3\frac{\partial \xi_2}{\partial X}.
$$
Then the Baker-Akhiezer function $\psi$ satisfies the following system
$$
\mathbb L\phi=\frac{\partial \psi}{\partial Y}, \quad \mathbb A\psi=\mu \psi.
$$
Compatibility condition for these two equations
$$
\big[-\frac{\partial}{\partial Y}+\mathbb L, \mathbb A\big]=0,
$$
 is equivalent to the system of nonlinear PDEs for $\u$ and $\w$:
$$
\frac{3}{4}\u_Y-\w_X=0,\quad \w_Y+\frac{1}{4}(\u_{XXX}+6\u\u_X)=0.
$$
By excluding $\w$ from the last two equations, one gets the following form of the Boussinesq equation \cite{Dub1981}:
\begin{equation}\label{eq:bous}
3\u_{YY}+(6\u\u_X+\u_{XXX})_X=0.
\end{equation}
The Boussinesq equation originated in the classical work of Boussinesq \cite{Bous} in 1872. It is a nonlinear integrable generalization of the wave equation, describing propagation of waves in weakly dispersive and weakly nonlinear media. The Lax pair for this equation was constructed by Zakharov \cite{Zakh} a century later. As indicated in \cite{Dub1981}, the function
\begin{equation}\label{eq:solbous}
\u(X, Y)=2\frac{\partial^2}{\partial X^2}\ln \theta(XU+YV+z_0)+c,
\end{equation}
solves the Boussinesq equation \eqref{eq:bous}, where $\theta(z)$ is the theta-function associated with $\surf$ and $U$ and $V$ are the $b$-periods of the differentials of the second kind \eqref{eq:UV} with a unique pole at $Q_0$, and $z_0, c\in\mathbb C$ are constants. As an application of Theorem \ref{thm_ode}, we get the following proposition.

\begin{proposition}\label{thm:bousappl} Consider a continuous family of solutions $\u(X, Y)$ \eqref{eq:solbous} of the Boussinesq equation \eqref{eq:bous}, constructed on the Riemann surfaces $\surf$ and an arbitrary selected point $Q_0\in\surf$ parameterized by the values of $x\in\X\subset \mathbb C\setminus\{0,1\},$ where $U$ and $V$ are given by \eqref{eq:UV}.

\begin{itemize}
\item[(i)] If the period $U$ (the wave number) is the same for all solutions $\u(y, t)$ in the family, then the projection $y_0=u(Q_0)$ of the point $Q_0$ on the $u$-sphere as a function of $x$ satisfies equation \eqref{ode} with  $n=0$.

\item[(ii)] If the period $V$ (the frequency) is the same for all solutions $\u(y, t)$ in the family, then the projection $y_0=u(Q_0)$ of the point $Q_0$ on the $u$-sphere as a function of $x$ satisfies equation \eqref{ode} with  $n=1$.
        \end{itemize}
\end{proposition}

{\it Proof.} This is a direct consequence of formula \eqref{eq:solbous} and Theorem \ref{thm_ode}.
$\Box$

Explicit forms of equation   \eqref{ode} for $n=0$ and $n=1$, used in (i) and (ii) of Proposition \ref{thm:bousappl} are given in Examples \ref{example_n0} and \ref{example_n1} respectively.
\\
Conversely, Theorem \ref{thm_unicity} implies the following proposition.
\begin{proposition}
\label{thm_bous_UV}
Let $x_0\in\mathbb C\setminus\{0,1\}$, and let $\mathcal L(x_0)$ be the surface corresponding to the elliptic curve given by \eqref{curve} with $x=x_0.$ Let $\{a,b\}$ be a canonical homology basis on $\mathcal L(x_0)$ such that the projections $u(a)$ and $u(b)$ do not intersect a certain neighbourhood $\hat\X$ of $x_0$.
Let $\Omega^{(n+1)}_{Q_0}(x_0)$ be the normalized meromorphic differential of the second kind defined by \eqref{eq:normsecond} with $n\geq 0$ on $\mathcal L(x_0)$ for some choice of a point $Q_0(x_0)\in\mathcal L(x_0)$ different from ramification points.  Let $\hat u$ be a solution to the Boussinesq equation defined by \eqref{eq:solbous} on $\mathcal L(x_0)$ corresponding to the point $Q_0\in\mathcal L(x_0)$ with $U$ and $V$ defined by \eqref{eq:UV}.
\begin{enumerate}
\item[(i)]
Assume that $n=0$ and that $Q_0(x_0)$ is such that $\omega^{(1)}(Q_0(x_0))\neq 0$ where $\omega^{(1)}$ is defined by \eqref{omega-Taylor}.
Then there exists a unique local continuous deformation of the solution $\hat u$ preserving the wave number $U$ \eqref{eq:UV}.
\item[(ii)]
Assume that $n=1$ and that $Q_0(x_0)$ is such that $\omega^{(2)}(Q_0(x_0))\neq 0$ where $\omega^{(2)}$ is defined by \eqref{omega-Taylor}. Then there exists a unique local continuous deformation of the solution $\hat u$ preserving the frequency $V$ \eqref{eq:UV}.
\end{enumerate}
These deformations correspond to the isoperiodic deformations $(\surf, \Omega^{(n+1)}_{Q_0}(x))$ of the pair \linebreak $(\mathcal L(x_0), \Omega^{(n+1)}_{Q_0}(x_0))$  for $x$ in a neighbourhood $\X\subset \hat\X$ of $x_0$.

\end{proposition}

A solution  $\u$ \eqref{eq:solbous}  is periodic in $X$ with period $T$ if and only if there exists a real number $T>0$, such that $TU$ is a lattice point of the Jacobian of $\surf$. Similarly, $\u$ is periodic in $Y$ with period $T$ if and only if there is a real number $T>0$, such that $TV$ is a lattice point of the Jacobian of $\surf$. The lattice of the Jacobian of $\surf$ is generated by $(1, \tau)$, see \eqref{eq:tau}.

Observe the following relationship between isoperiodic deformations differentials of the second kind on Riemann surfaces and periodic solutions of the Boussinesq equation.  Let $\u(X, Y)$ be a solution  \eqref{eq:solbous} of the Boussinesq equation \eqref{eq:bous}, constructed from the Riemann surface $\mathcal L(x_0)$ and its selected point $Q_0$ for which $U$ \eqref{eq:UV} is real. Then, if $\surf,$ $x\in\X\subset \mathbb C\setminus\{0,1\}$ where $x_0\in \X$ is a family of surfaces from Theorem \ref{thm_ode} with respect to the differential $\Omega=\Omega_{Q_0}^{(1)}$ \eqref{eq:normsecond}, then solutions $\u$ of the Boussinesq equation constructed from surfaces $\surf$ form a continuous family of solutions periodic in $X$ with period $T=1/U$.
\begin{proposition}
\label{prop_BT2}
Consider a continuous family of solutions $\u(X, Y)$ \eqref{eq:solbous} of the Boussinesq equation \eqref{eq:bous}, constructed from a family of Riemann surfaces $\mathcal L(x)$ parameterized by $x\in\X\subset \mathbb C\setminus\{0,1\}$,  with $\X$ being  connected, and their selected points $Q_0=Q_0(x)\in\surf$.
\begin{itemize}
\item[(i)] Assume that each solution $\u(X, Y)$ from the family is periodic in $X$ with a period $T\in\mathbb R$. Assume also that  there exists $x_0\in\X$ such that $U(x_0)$ \eqref{eq:UV} is real for the surface  $\mathcal L(x_0)$. Then $U$ \eqref{eq:UV} is real for all surfaces $\mathcal L(x)$ with $x\in\X\subset \mathbb C\setminus\{0,1\}$ and the projection $y_0$ of $Q_0$ on the $u$-sphere as a function of $x$ satisfies equation \eqref{ode} with $n=0$.

    \item[(ii)] Assume that each solution $\u(X, Y)$ from the family is periodic in $Y$ with a period $T\in\mathbb R$. Assume also that  there exists $x_0\in\X$ such that $V(x_0)$ \eqref{eq:UV} is real for the surface  $\mathcal L(x_0)$. Then $V$ \eqref{eq:UV} is real for all surfaces $\mathcal L(x)$ with $x\in\X\subset \mathbb C\setminus\{0,1\}$  and  the projection $y_0$ of $Q_0$ on the $u$-sphere as a function of $x$ satisfies equation \eqref{ode} with $n=1$.
        \end{itemize}
\end{proposition}
{\it Proof.} If solutions $\u(X,Y)$ are periodic  in $X$ with period $T\in\mathbb R$, then $TU$ is a lattice vector for all $x\in\mathcal X$. If in addition  $U(x_0)$ is real, then $TU(x_0)$ is an integer and thus, due to the continuity of $U(x)$ and connectedness of $\X$,  the quantity $TU(x)$ is integer for all $x\in\X$ and  $TU(x)=TU(x_0)$. Thus $U(x)$ is constant for all curves in the family $\surf$. Then, Theorem \ref{thm_ode} is applicable to the family $\surf$ with $n=0$ and thus $y_0$ satisfies \eqref{ode} with $n=0.$ The proof of the second part of the theorem is analogous with $n=1$.
$\Box$

Conversely, from Theorem \ref{thm_unicity}, we obtain the following proposition.
\begin{proposition}
\label{thm-periodic-deformation}
Let $x_0\in\mathbb C\setminus\{0,1\}$, and let $\mathcal L(x_0)$ be the surface corresponding to the elliptic curve given by \eqref{curve} with $x=x_0.$  Let $\{a,b\}$ be a canonical homology basis on $\mathcal L(x_0)$ such that the projections $u(a)$ and $u(b)$ do not intersect a certain neighbourhood $\hat\X$ of $x_0$.
Let $\Omega^{(n+1)}_{Q_0}(x_0)$ be the normalized meromorphic differential of the second kind defined by \eqref{eq:normsecond} with $n\geq 0$ on $\mathcal L(x_0)$ for some choice of a point $Q_0(x_0)\in\mathcal L(x_0)$ different from ramification points. Let $\hat u$ be a solution to the Boussinesq equation defined by \eqref{eq:solbous} on $\mathcal L(x_0)$ corresponding to the point $Q_0\in\mathcal L(x_0)$.
\begin{enumerate}
\item[(i)]
Assume that $n=0$ and that the surface and the point $Q_0(x_0)\in\mathcal L(x_0)$ are such that $U$ \eqref{eq:UV} is real and $\omega^{(1)}(Q_0(x_0))\neq 0$ where $\omega^{(1)}$ is defined by \eqref{omega-Taylor}. In this case, the solution $\hat u$ is periodic in $X$ with period $1/U$ and there is a unique continuous deformation  $\hat u(x, X,Y)$  of $\hat u(X,Y)$ such that for each $x$ in some heighbourhood $\X\subset\hat\X$ of $x_0$ the function $\hat u(x,X,Y)$  is a periodic in $X$ solution to the Boussinesq equation \eqref{eq:bous} with the same period $1/U$.
\item[(ii)]
Assume that $n=1$ and that the surface and the point $Q_0(x_0)\in\mathcal L(x_0)$ are such that $V$ \eqref{eq:UV} is real and $\omega^{(2)}(Q_0(x_0))\neq 0$ where $\omega^{(2)}$ is defined by \eqref{omega-Taylor}. In this case, the solution $\hat u$ is periodic in $Y$ with period $1/V$ and there is a unique continuous deformation  $\hat u(x, X,Y)$  of $\hat u(X,Y)$ such that for each $x$ in some neighbourhood $\X\subset\hat\X$ of $x_0$ the function $\hat u(x,X,Y)$  is a periodic in $Y$ solution to the Boussinesq equation \eqref{eq:bous} with the same period $1/V$.
\end{enumerate}
\end{proposition}
{\it Proof.} Given that $U$ is real for $x=x_0$, Theorem \ref{thm_unicity} gives a unique deformation of the pair $(\mathcal L(x_0), \Omega^{(n+1)}_{Q_0})$ preserving the value of $U$. Thus $U$ stays real for any curve from the isoperiodic family $(\surf, \Omega^{(n+1)}_{Q_0})$. This implies that solutions $\hat u$ \eqref{eq:solbous} constructed from $\surf$ for $x\in\X$ is also periodic in $X$ and the period is $1/U.$
$\Box$

\smallskip

\begin{remark}\label{rem:real} According to Remark \ref{rem:realell}, assuming that the homology cycles $a$ and $b$ are selected as in that remark, it is enough to select real $x_0$ and $y_0(x_0)$ such that $y_0(y_0-1)(y_0-x)>0$ in order for $U(x_0)$ and $V(x_0)$ to be real, that is in order for the assumptions of Propositions \ref{prop_BT2} and \ref{thm-periodic-deformation} concerning the reality of $U(x_0)$ of $V(x_0)$ to hold.
\end{remark}

\smallskip

\begin{remark}\label{rem:eff}
By applying Dubrovin's effectivization procedure from \cite{Dub1981}, Ch. IV (see formulas (4.1.13') and (4.2.5)) we get the relationship between $U$ and $V$:
$$
V=\pm \frac{4\pi\i\sqrt{3}}{3}U^2\sqrt{\frac{\theta_1^{'''}(0)}{\theta_1'(0)}}.
$$
\end{remark}

\bigskip
\bigskip

{\bf Acknowledgements.}   V. D. acknowledges with gratitude the support from the Simons Foundation grant no. 854861 and the Science Fund of Serbia grant IntegraRS. V. S. gratefully acknowledges support from the Natural Sciences and Engineering Research Council of Canada through a Discovery grant and from the University of Sherbrooke.
\\

%{\bf Data availability statement.}
%The authors declare that there is no data associated with the manuscript.
%


\begin{thebibliography}{99}

\addcontentsline{toc}{section}{Bibliography}

\bibitem{parts} G. E. Andrews, {\it The theory of partitions.} Cambridge Math. Lib.
Cambridge University Press, Cambridge, 1998, xvi+255 pp.

\bibitem{Bous} J. Boussinesq, \textit{
Th\'eorie des ondes et des remous qui se propagent le long d'un canal
rectangulaire horizontal, en communiquant au liquide contenu dans
ce canal des vitesses sensiblement pareilles de la surface au fond},
Journal de math\'ematiques pures et appliqu\'ees 2e s\'erie, tome 17 (1872), p. 55--108.



\bibitem{isoharmonic} V. Dragovi\'c, V. Shramchenko, \textit{Isoharmonic deformations and constrained Schlesinger systems} Advances in Mathematics, Vol. 499, July 2026. 

\bibitem{DS2019} V. Dragovi\'{c}, V. Shramchenko, \textit{Algebro-geometric approach to an Okamoto transformation, the Painlev\'e VI and Schlesinger equations}, Annales Henri Poincar\'e, 2019,  20, no. 4, 1121--1148.

\bibitem{Dub1981} B. A. Dubrovin, \textit{Theta functions and non-linear equations}, Russ. Math. Surv. 1981,  36, No. 2 11--92.


\bibitem{DubNat} B. A. Dubrovin, S. M.  Natanzon, \textit{Real two-zone solutions of the sine-Gordon equation}, Funktsional. Anal. i Prilozhen., (1982), 16:1,  27--43.


\bibitem{Fay92} J. Fay, {\it Kernel functions, analytic torsion, and moduli spaces}, Memoirs of the AMS \textbf{96} (1992), no 464.

\bibitem{GS} P. G. Grinevich, M. U. Schmidt, {\it Period preserving nonisospectral flows and the moduli space of periodic solutions of soliton equations: The nonlinear Schr\"odinger equation}, Physica D, 87:73--98, 1995.

\bibitem{Ince} E. L. Ince, {\it Ordinary differential equations} Dover, New York (1956).



\bibitem{KokoKoro2} A. Kokotov, D. Korotkin, {\it Isomonodromic tau-function of Hurwitz Frobenius manifolds and its applications.}
Int. Math. Res. Not. 2006, Art. ID 18746, 34 pp.

\bibitem{KokoKoro} A. Kokotov, D. Korotkin, {\it A new hierarchy of integrable systems associated to Hurwitz spaces}, Philos. Trans. R. Soc. Lond. Ser. A Math. Phys. Eng. Sci. {\bf 366} (2008), no. 1867, 1055--1088.


\bibitem{Picard} E. Picard, {\it  M\'emoire sur la th\'eorie des
fonctions alg\'ebriques de deux variables}, Journal de Liouville
{\bf 5} 1889, 135--319.

\bibitem{Rauch} Rauch, H.E., Weierstrass points, branch points, and moduli of Riemann surfaces, Comm. Pure
Appl. Math. {\bf 12} 543-560 (1959).

\bibitem{Sp1957} G. Springer, {\it Introduction to Riemann Surfaces}, AMS Chelsea Publishing, 1957, p. 309.

\bibitem{Zakh} V. E. Zakharov, \textit {On stochastization of one-dimensional chains of
nonlinear oscillators}, Sov. Phys.--JETP, Vol. 38, No. 1, January 1974, 108--110.

\end{thebibliography}
\end{document}